\documentclass[preprint,9pt]{elsarticle}

\usepackage{lineno,hyperref}
\usepackage{amsmath}
\usepackage{amsfonts}
\usepackage{latexsym}
\usepackage{amsthm}
\usepackage{amssymb, ifthen, enumerate}
\usepackage{color}
\usepackage{ulem}
\usepackage{longtable}
\usepackage{listings}
\usepackage{xcolor}
\modulolinenumbers[5]
\definecolor{codegray}{rgb}{0.6,0.6,0.6}
\definecolor{codenormal}{rgb}{0.35,0.35,0.35}
\definecolor{codeblack}{rgb}{0,0,0}
\definecolor{backcolour}{RGB}{250,250,250}

\lstdefinestyle{gapstyle}{
  backgroundcolor=\color{backcolour}, 
  commentstyle=\itshape\color{codegray},
  keywordstyle=\color{codenormal},
  numberstyle=\tiny\color{codegray},
  basicstyle=\ttfamily\footnotesize,
  breakatwhitespace=false, 
  breaklines=true,                 
  captionpos=b,                    
  keepspaces=true,                 
  numbers=none,                    
  numbersep=5pt,                  
  showspaces=false,                
  showstringspaces=false,
  showtabs=false,                  
  tabsize=2
}

\journal{Journal of Algebra, Computational Section}

\begin{document}

\begin{center}
\Large{\textbf{Finite simple groups acting with fixity $4$ and their occurrence as groups of automorphisms of Riemann surfaces (extended version)}}

\vspace{0.2cm} \small{Patrick Salfeld and Rebecca Waldecker}\footnote{MLU Halle Wittenberg, Institut f\"ur Mathematik, 06099 Halle, Germany, rebecca.waldecker@mathematik.uni-halle.de}
\end{center}

\normalsize
\vspace{0.5cm}
\noindent
\textbf{Abstract.}

In previous work, all finite simple groups that act with fixity 4 have been classified.
In this article we investigate which ones of these groups act faithfully on a compact Riemann surface of genus at least $2$ with fixity four in total and in such a way that
fixity 4 is exhibited on at least one orbit. This is an extended version of the submitted article, including our \texttt{GAP} code.


\newtheorem{lemma}{Lemma}[section]
\newtheorem{thm}[lemma]{Theorem}
\newtheorem{prop}[lemma]{Proposition}
\newtheorem{cor}[lemma]{Corollary}
\newtheorem{rem}[lemma]{Remark}
\newtheorem{hyp}[lemma]{Hypothesis}
\newtheorem{definition}[lemma]{Definition}
\newtheorem{ex}[lemma]{Example}
\newtheorem{Notation}[lemma]{Notation}

\renewcommand{\qed}{\hspace*{\fill}\rule{1.5ex}{1.5ex}}
\newcommand{\N}{\mathbb{N}}	
\newcommand{\R}{\mathbbm{R}}	
\newcommand{\Q}{\mathbbm{Q}}	
\newcommand{\Z}{\mathbbm{Z}}	
\newcommand{\F}{\mathbbm{F}}	
\newcommand{\C}{\mathbbm{C}}	
\newcommand{\Ri}{\mathcal{R}}     
\newcommand{\Alt}{\textrm{Alt}}     
\newcommand{\Sym}{\textrm{Sym}}     
\newcommand{\mc}[1]{\mathcal{#1}}	
\newcommand{\mf}[1]{\mathfrak{#1}}	
\newcommand{\ug}{\leqslant}		
\newcommand{\eug}{\lneqq}			
\newcommand{\gu}{\geqslant}		
\newcommand{\gue}{\gneqq}			
\newcommand{\marg}[1]{\marginpar{\fbox{\centering\footnotesize\texttt{#1}}}}
\newcommand{\triv}{\mathbbm 1}

\newcommand{\Menge}[2]{\left\{ #1 \;\middle|\; #2 \right\}}
\newcommand{\Erzeug}[2]{\left\langle #1 \;\middle|\; #2 \right\rangle}
\newcommand{\Aut}{\textrm{Aut}}
\newcommand{\Kern}{\textrm{Kern}}
\newcommand{\Bild}{\textrm{Im}}
\newcommand{\Fix}{\textrm{Fix}}
\newcommand{\M}{\textrm{M}}
\newcommand{\id}{\textrm{id}}
\newcommand{\GL}{\textrm{GL}}
\newcommand{\SL}{\textrm{SL}}
\newcommand{\Hom}{\textrm{Hom}}
\newcommand{\End}{\textrm{End}}
\newcommand{\Hol}{\textrm{Hol}}
\newcommand{\Mer}{\textrm{Mer}}
\newcommand{\Rang}{\textrm{Rang}}
\newcommand{\eps}{\varepsilon}
\renewcommand{\tilde}{\widetilde}
\newcommand{\mult}{mult}
\newcommand{\Gal}{\textrm{Gal}}
\newcommand{\PSL}{\textrm{PSL}}
\newcommand{\PSp}{\textrm{PSp}}
\newcommand{\POm}{\textrm{P$\Omega$}}
\newcommand{\PGL}{\textrm{PGL}}
\newcommand{\PSU}{\textrm{PSU}}
\newcommand{\PGU}{\textrm{PGU}}
\newcommand{\GU}{\textrm{GU}}
\newcommand{\GO}{\textrm{GO}}
\newcommand{\SU}{\textrm{SU}}
\newcommand{\Sp}{\textrm{Sp}}
\newcommand{\GF}{\textrm{GF}}
\newcommand{\Sz}{\textrm{Sz}}
\newcommand{\J}{\textrm{J}}
\newcommand{\D}{\textrm{D}}
\newcommand{\G}{\textrm{G}}
\newcommand{\Syl}{\textrm{Syl}}
\newcommand{\ggT}{\textrm{gcd}}
\newcommand{\kgV}{\textrm{kgV}}
\newcommand{\Irr}{\textrm{Irr}}
\newcommand{\FO}{\textrm{fix}_{\Omega}}
\newcommand{\FD}{\textrm{fix}_{\Delta}}
\makeatletter 

\newcommand{\be}[1]{\textbf{#1}}
\renewcommand{\labelenumi}{(\textit{\alph{enumi}})}
\newcommand{\abs}[1]{\left| #1 \right|}
\newcommand{\SSD}[1][G, \mf g,\mf g_0]{\left[ #1 \;\bigr |\; [m_1,n_1],\ldots,[m_r,n_r]\right]}
\newcommand{\ssd}{\SSD}
\newcommand{\SD}{\SSD[\mf g_0]}
\newcommand{\sd}{\SD}
\newcommand{\SDo}[1][G, \mf g,\mf g_0]{\left[ #1 \;\bigr |\; 0\right]}
\newcommand{\sdo}{\SDo}
\newcommand{\ld}{\LD}
\newcommand{\cog}{{\mf g_0}}

\newcommand{\prf}[1]{
\noindent Proof.
#1
\qed\\[0.5eM]
}

\section{Introduction}

Suppose that a group $G$ acts faithfully and transitively on a finite set $\Omega$.
Following Ronse (see \cite{Ronse}), we say that $G$ acts with fixity $k \in \N$ on $\Omega$
  if and only if $k$ is the maximum number of fixed points of non-trivial elements from $G$.
While the characterisation of finite permutation groups that act with low fixity
is interesting in its own right, we had applications in mind that are motivated from
the theory of Riemann surfaces. This is why ``low fixity'' means fixity 2, 3 or 4. This is also why we do not stop at the
classification of finite simple groups that act with fixity 2, 3 or 4, but we always
go one step further and investigate which ones of our examples actually occur as subgroups of the automorphism group of a compact Riemann surface of genus at least $2$.
We will not go into details here, but instead refer to
Magaard and V\"olklein (\cite{MV}) and to our previous articles \cite{SW} and \cite{SW3} for more background.

Our strategy in this article is similar to \cite{SW3},
which means that we build on the classification of the finite simple groups that act with fixity 4,
then we only consider the cases where point stabilisers are cyclic and finally derive the branching data for possible curves. With the applications in mind we suppose, for the branching data, that there is at least one non-regular orbit on which the group acts with fixity 4.
Otherwise this possibility would have been covered in earlier work. 
Some groups allow for actions with fixity 2, 3, or 4, and then there might be several non-regular orbits and the group acts with fixity 2, 3 or 4 on them, respectively. We will refer to this behaviour as \textbf{``mixed fixity action''} and discuss it in more detail later.

First we state
the main theorem:

\begin{thm}\label{main4fpH}
Suppose that $G$ is a finite simple group that acts as a group of automorphisms on a compact Riemann surface
$X$ of genus at least $2$.
Suppose further that there is at least one orbit on which $G$ acts with fixity 4 and that $G$ acts with fixity at most 4 in total on $X$.
Then $G$ acts with one of the branching data given in Table \ref{DataList}.

Conversely, if $l$ is a list from Table \ref{DataList} and $l$ is a Hurwitz datum, then there exists a compact Riemann surface $X$ of genus $g$ at least 2 on which the group $G$ in $l$ acts faithfully, and for all choices of $X$, it is true that $G$ acts with fixity 4 on at least one orbit of $X$ and with fixity at most 4 in total.


\end{thm}


In the following table, $q$ always denotes a prime power. We also remark that, in particular, $g_0 \ge 0$ if a Hurwitz datum has at least three branchings, and otherwise $g_0 \ge 1$. But there are some exceptions which we mention in the "Remark" column. 

In Lines 8a--8i we set $\alpha:=q^2+\sqrt{2 \cdot q}+1$ and $\beta:=q^2-\sqrt{2 \cdot q}+1$, and $q$ is a power of $2$. The column "Fixity" indicates the fixity with which $G$ acts on the non-regular orbits, in the order in which they appear in the Hurwitz datum. 

\begin{center}
\begin{longtable}{llll}
	
	Line&Hurwitz datum&Remark&Fixity\\[0.4ex]\hline\hline
	\\[-2ex]

        1a&$[\Alt_7,g,g_0 \mid [5,1]]$& & 4\\
        1b&$[\Alt_7,g,g_0 \mid [5,1],[7,1] ]$& &4,\,3\\
        \\[-2ex]\hline\\[-2ex]
		2a&$[\PSL_2(7),g,g_0 \mid [2,1]]$&$g_0 \ge 2$&4\\
		2b&$[\PSL_2(7),g,g_0 \mid [2,1],[3,1]]$& &4,\,2\\
		2c&$[\PSL_2(7),g,g_0 \mid [2,1],[3,2]]$&$g_0 \ge 1$&4,\,2\\
		2d&$[\PSL_2(7),g,g_0 \mid [2,1],[7,1]]$& &4,\,3\\
		2e&$[\PSL_2(7),g,g_0 \mid [2,1],[3,1],[7,1]]$& &4,\,3,\,2\\
		2f&$[\PSL_2(7),g,g_0 \mid [2,1],[3,2],[7,1]]$& &4,\,3,\,2\\
          \\[-2ex]\hline\\[-2ex]
	        3a&$[\PSL_2(8),g,g_0 \mid [2,1]]$&$g_0 \ge 2$&4\\
		3b&$[\PSL_2(8),g,g_0 \mid [2,1],[7,1]]$& &4,\,2\\
		3c&$[\PSL_2(8),g,g_0 \mid [2,1],[7,2]]$& &4,\,2\\
		3d&$[\PSL_2(8),g,g_0 \mid [2,1],[9,1]]$& &4,\,2\\
		3e&$[\PSL_2(8),g,g_0 \mid [2,1],[9,2]]$& &4,\,2\\
		3f&$[\PSL_2(8),g,g_0 \mid [2,1],[7,1],[9,1]]$& &4,\,2,\,2\\
		3g&$[\PSL_2(8),g,g_0 \mid [2,1],[7,2],[9,1]]$& &4,\,2,\,2\\
		3h&$[\PSL_2(8),g,g_0 \mid [2,1],[7,1],[9,2]]$& &4,\,2,\,2\\
		3i&$[\PSL_2(8),g,g_0 \mid [2,1],[7,2],[9,2]]$& &4,\,2,\,2\\
  \\[-2ex]\hline\\[-2ex]
                4a&$[\PSU_4(3),g,g_0 \mid [5,1]]$& &4\\

                4b&$[\PSU_4(3),g,g_0 \mid [5,1],[7,1]]$& &4,\,3\\
                \\[-2ex]\hline\\[-2ex]
5a&$[\PSL_2(q),g,g_0 \mid [\frac{q-1}{4},1]]$& $q \ge 9$, $q \equiv 1 \mod 4$&4\\
5b&$[\PSL_2(q),g,g_0 \mid [\frac{q+1}{4},1]]$& $q \ge 9$, $q \equiv -1 \mod 4$&4\\
5c&$[\PSL_2(q),g,g_0 \mid [\frac{q-1}{4},1],[\frac{q+1}{2},1]]$& $q \ge 9$, $q \equiv 1 \mod 4$&4,\,2\\
5d&$[\PSL_2(q),g,g_0 \mid [\frac{q-1}{4},1],[\frac{q+1}{2},2]]$& $q \ge 9$, $q \equiv 1 \mod 4$&4,\,2\\
5e&$[\PSL_2(q),g,g_0 \mid [\frac{q+1}{4},1],[\frac{q-1}{2},1]]$& $q \ge 9$, $q \equiv -1 \mod 4$&4,\,2\\
5f&$[\PSL_2(q),g,g_0 \mid [\frac{q+1}{4},1],[\frac{q-1}{2},2]]$& $q \ge 9$, $q \equiv -1 \mod 4$&4,\,2\\
\\[-2ex]\hline\\[-2ex]
6&$[\PSp_4(q),g,g_0 \mid [\frac{q^2+1}{(2,q^2-1)},1]]$& $q \ge 3$&4\\
\\[-2ex]\hline\\[-2ex]
7&$[\POm^-_8(q),g,g_0 \mid [\frac{q^4+1}{(2,q^4-1)},1]]$& &4 \\
\\[-2ex]\hline\\[-2ex]
	        8a&$[\Sz(q),g,g_0 \mid [\alpha,1]]$& $q \ge 8$&4\\
                8b&$[\Sz(q),g,g_0 \mid [\beta,1]]$&$q \ge 8$&4\\
	        8c&$[\Sz(q),g,g_0 \mid [\alpha,1],[\beta,1]]$& $q \ge 8$&4,\,4\\
		8d&$[\Sz(q),g,g_0 \mid [\alpha,1],[q-1,1]]$& $q \ge 8$&4,\,2\\
	        8e&$[\Sz(q),g,g_0 \mid [\alpha,1],[q-1,2]]$& $q \ge 8$&4,\,2\\
	        8f&$[\Sz(q),g,g_0 \mid [\beta,1],[q-1,1]]$& $q \ge 8$&4,\,2\\
	        8g&$[\Sz(q),g,g_0 \mid [\beta,1],[q-1,2]]$& $q \ge 8$&4,\,2\\
	        8h&$[\Sz(q),g,g_0 \mid [\alpha,1],[q-1,1],[\beta,1]]$& $q \ge 8$&4,\,2,\,4\\
	        8i&$[\Sz(q),g,g_0 \mid [\alpha,1],[q-1,2],[\beta,1]]$& $q \ge 8$&4,\,2,\,4\\
         \\[-2ex]\hline\\[-2ex]
		9&$[{^3}\D_4(q),g,g_0 \mid [q^4-q^2+1,1]]$& &4\\
  \\[-2ex]\hline\\[-2ex]
		10&$[{^2}\G_2(q),g,g_0 \mid [\frac{q-1}{2},1]]$& $q \ge 27$&4\\
  \\[-2ex]\hline\\[-2ex]

		11&$[\M_{11},g,g_0 \mid [5,1]]$& &4\\
  \\[-2ex]\hline\\[-2ex]

		12a&$[\M_{22},g,g_0\mid[5,1]]$& &4\\
		12b&$[\M_{22},g,g_0\mid[5,1],[7,1]]$& &4,\,3\\
  \\[-2ex]\hline\\[-2ex]
	    13&$[\J_1,g,g_0\mid[15,1]]$& &4\\

		\\[-2ex]\hline\hline\\[-2ex]

	\caption{\small Hurwitz data for group actions with only one non-regular orbit and fixity 4 or with up to five non-regular orbits and mixed fixity, $g \ge 2$.}
	\label{DataList}
	
\end{longtable}
\end{center}


We give an example that illustrates how to use Table \ref{DataList}.

In Lines 4a and 4b we see the group $\PSU_4(3)$, and we first look at Line 4a with the list $[\PSU_4(3),g,g_0 \mid [5,1]]$. This is a potential branching datum for the group $\PSU_4(3)$, acting as a group of automorphisms on a compact Riemann surface with genus $g$ and cogenus $g_0$ and such that there is exactly one non-regular orbit.
The "Fixity" column indicates that $G$ acts with fixity $4$ on this orbit,
and the Hurwitz datum indicates that the point stabilisers have order $5$ and fix exactly four points, all in this unique non-regular orbit.

In Line 4b with the list $[\PSU_4(3),g,g_0 \mid [5,1],[7,1]]$, we have the same group, this time with what we call a \textbf{mixed fixity action}:
The Hurwitz datum shows us two non-regular orbits, and point stabilisers of order $5$ and order $7$, respectively. Together with the "Fixity" column, we see that the first non-regular orbit belongs to a fixity $4$ action with point stabilisers of order $5$, and the second non-regular orbit belongs to a fixity $3$ action with point stabilisers of order $7$.\\


For the proof of our main theorem,
we recall some relevant definitions, and then we refer to Theorem 1.3 of \cite{BHMSW} for the specific groups to consider. In the analysis
we follow the strategy developed in \cite{SW}: First, in Section 3, we analyse series of groups and find generating sets with specific properties. Then we explain the arguments for the proof of Theorem \ref{main4fpH} in Section 4. Some individual groups can be treated with \texttt{GAP} (\cite{GAP4}), and we show and comment on the code at the end of this article.\\  


\textbf{Acknowledgements.} We remember Kay Magaard with much appreciation and gratitude for
initiating the project on groups acting with low fixity and for
suggesting questions along the lines of this article.
We are also grateful to the referee because one of their comments gave us the idea to add information to Table 1. 


\section{Preliminaries}

Most of our notation is standard, we refer to \cite{Miranda1997} for background information on Riemann surfaces
and we only mention the following for clarity: \\
If $G$ is a group and $g,h \in G$, then we write \textbf{conjugation} from the right, so $g^h:=h^{-1}  g h$, and we write \textbf{commutators} as $[g,h]:=g^{-1}  h^{-1}  g  h$.\\
Throughout, we suppose that $X$ is a compact Riemann surface of genus $g\ge 2$ and that $G\le\Aut(X)$. Let $X/G=\Menge{ x^G}{ x\in  X}$ denote the \textbf{space of $G$-orbits on $X$} and let $ g_0\ge 1$ denote the \textbf{co-genus}, i.e. the genus of $X/G$.
Then all elements in $X/G$ are finite and their cardinalities (which we refer to as \textbf{orbit lengths}) divide $|G|$.
We recall that an orbit is \textbf{non-regular }if and only if some element of $G^\#$ fixes a point in it, which happens if and only if its length is strictly less than $|G|$. For all $x\in X$ we denote the point stabiliser of $x$ in $G$ by $G_x$.
In \cite{SW} we give more details about the notation. Here we just recall:
	\begin{itemize}
                \item $G$ has infinitely many regular orbits on $X$ and only finitely many non-regular orbits. (Proposition III.3.2, Theorem III 3.4 and
Corollary II.1.34 in \cite{Miranda1997}.)
		\item  $\abs {\Aut( X)}$ is finite and for all $x \in X$, the subgroup $G_x$ is cyclic. (Theorem VII.4.18 and Proposition
III.3.1 in \cite{Miranda1997}.)
	\end{itemize}
In particular, all groups considered here are finite.

\begin{definition}\label{einf1-not:stlit}

Given a finite group $G$, non-negative integers $g,g_0,r,m_1,n_1,...,$\\
$m_r,n_r$ and a list
$l:=[G, g, g_0 \mid [m_1,n_1],\ldots, [m_r,n_r]]$,
we refer to $l$ as a \textbf{Hurwitz datum} if and only if the Hurwitz formula
is satisfied:
$$2(g-1)=\abs G\left( 2(g_0-1)+\sum_{j=1}^r n_j\left(1-\frac 1{m_j}\right)\right).$$

We say that $G$ \textbf{acts with branching datum} $l$ on a Riemann surface $X$ if and only if
$X$ has genus $g$, $X/G$ has genus $g_0$ and for each $i \in \{1,...,r\}$ there are exactly $n_i$ non-regular orbits
of $X/G$
on which $G$ acts with point stabilisers of order $m_i$.
In such a case we refer to $X$ as a \textbf{witness for $l$}.
\end{definition}


Our starting point is the table from Theorem 1.3 in \cite{BHMSW}, we list
only the possibilities with cyclic point stabilisers, and we make the list more compact than in \cite{BHMSW}.
For convenience, we have a separate column that indicates mixed fixity, based on 
Theorem 1.2 and Lemmas 3.11 and 3.13 in \cite{MW} and Theorem 1.1 in \cite{MW3}. 
The number in brackets in the right column indicates the point stabiliser order for fixity 2 or 3 in the mixed cases. 

We remark that $q$ always denotes a prime power and that $q \ge 9$ in the generic $\PSL_2(q)$ case.

\begin{center}
\begin{longtable}{lll}

\mbox{}
	Group \(G\)                           & Possible point stabiliser order  & mixed fixity\\ [0.4ex]\hline\hline \\

	\(\Alt_7\)                            & \(5\) & 3 (7)\\

	\(\PSL_2(7) \cong \PSL_3(2)\)         & \(2\) & 2 (3) and 3 (7) \\
	
	\(\PSL_2(8)\)                         & \(2\)& 2 (7 or 9) \\

		\(\PSU_4(3)\)          & \(5\)       & 3 (7)                     \\[0.5ex]
		
		\(\PSL_2(q), q \equiv 1\) mod 4                      &  $\frac{q-1}{4}$ & 2 ($\frac{q+1}{2}$) \\[1ex]
	
	       \(\PSL_2(q), q \equiv -1\) mod 4                       &  $\frac{q+1}{4}$  & 2 ($\frac{q-1}{2}$) \\[1ex]

			\(\PSp_4(q)\)          &  $\frac{q^2+1}{(2,q^2-1)}$           &  none              \\[1ex]
		
		\(\Sz(q)\)                            &  $q+\sqrt{2q}+1$ or
		$q-\sqrt{2q}+1$              & 2 ($q-1$)           \\[0.5ex]
			\(\POm^-_8(q)\)          &    $\frac{q^4+1}{(2,q^4-1)}$ &  none   \\[1ex]
			\(^3D_4(q)\)          &   $q^4-q^2+1$          &   none         \\[0.5ex]
			\(^2G_2(q)\)          &  $\frac{q-1}{2}$   &  none                  \\[1ex]
	
	\(M_{11}\)                            & \(5\) & none \\

	\(M_{22}\)                            & \(5\)     & 3 (7)                                    \\

	\(J_1\)                               & \(15\)    &  none                               \\

		\\[-2ex]\hline\hline\\[-2ex]
	\caption{\small Finite simple groups that act with fixity 4 and cyclic point stabilisers.}
	\label{GroupList}

\end{longtable}
\end{center}



\section{Analysis of series of groups}

In order to prove Theorem \ref{main4fpH}, we need to prove that every group satisfying the first hypothesis appears in Table \ref{DataList}, and only with the branching data that is given there.
This requires us to look at every group (or series of groups) individually. 
Conversely, we need to look at each group (or series of groups) from the table again and prove that there exists a witness 
in the sense of the previous section.
Often, the necessary calculations can be done with \texttt{GAP} (see \cite{GAP4}) and we will show the code that we used, going through the groups in the order in which they appear in Table \ref{GroupList}.

Here we look at generic examples, which means series of groups, and for each group we prove an individual lemma on our way to Theorem \ref{main4fpH}. Some preparation is needed, in particular we will use several times that the existence of a witness for a Hurwitz datum is equivalent to the existence of a generating set for the corresponding group $G$ with specific properties (originally from \cite{Broughton1990}):

\begin{lemma}\label{crit}
Suppose that $l=[G, g, g_0 \mid [m_1,n_1],\ldots, [m_r,n_r]]$ is a Hurwitz datum.
Then $G$ has a witness $X$ with respect to $l$ if and only if there exist
elements
 $a_1,\ldots,a_{g_0},b_1,\ldots,b_{g_0},c_{1,1},\ldots,c_{1,n_1},\ldots,c_{r,1},
\ldots,c_{r,n_r}\in G$ satisfying the following conditions:

\begin{enumerate}
\item
For all $j\in\{1,\ldots,r\}$ and all $i\in\{1,\ldots,n_j\}$ it is true that $o(c_{j,i})=m_j$,
    \item
    $\prod_{k=1}^{g_0}[a_k,b_k]\cdot \prod_{i=1}^{n_1}c_{1,i}\cdots \prod_{i=1}^{n_r}c_{r,i}=1$, and
    \item
    $\langle a_1,\ldots, a_{g_0},b_1,\ldots,b_{g_0},c_{1,1},\ldots,c_{r,n_r}\rangle=G$. \end{enumerate}

Referring to (c) we often denote the set as a generating tuple for $G$, namely as $((a_1,...,a_{g_0}),(b_1,...,b_{g_0}),(c_{1,1},...,c_{r,n_r}))$, and in the special case of cogenus $0$ we write 
$((~),(~),(c_{1,1},...,c_{r,n_r}))$.

Furthermore, if $G$ has a witness $X$ with respect to $l$ and if $h\in G^\#$ fixes a point in $X$, then $h$ is conjugate to a power of one of $c_{1,1},\ldots,c_{r,n_r}$ in $G$.
\end{lemma}

If $C_1, C_2, C_3$ are conjugacy classes of elements of a group $G$ of orders $m_1$, $m_2$ and $m_3$, respectively, then we define, for all  $i,j,l \in \{1,2,3\}$, coefficients 
$d_{ijl} := |\{(g_i, g_j) \in C_i \times C_j \mid g_i \cdot g_j = g_l \in C_l\}|$.
We refer to Remark 3.4 in \cite{SW} for more details.

\begin{lemma}\label{max}
Suppose that $n_1,n_2 \in \{0,1\}$, $1 \le n_1 + n_2 \le 2$ and that $l := [G,g,1 \mid [m_1,n_1], [m_2,n_2]]$ is a Hurwitz datum. 
Suppose further that $\bar{n_1},\bar{n_2} \in \{1,2\}$ are such that $\bar{n_1} +\bar{n_2} = 3$ and that $\bar{l} := [G,\bar{g},0 \mid [m_1,\bar{n_1}], [m_2,\bar{n_2}]]$ is a Hurwitz datum. 
Let $C_1,C_2$ denote a conjugacy class of $G$ of elements of order $m_1,m_2$, respectively, and suppose that $d_{212} \neq 0 \neq d_{121}$.
Then the properties (a) and (b) of Lemma \ref{crit} are satisfied for $l$ and $\bar{l}$.

If, in addition, every maximal subgroup $M$ of $G$ intersects $C_1$ or $C_2$ trivially (or equivalently if $|M|$ is never divisible by the lowest common multiple of $m_1$ and $m_2$), then $l$ and $\bar{l}$ have witnesses. 

\end{lemma}

\begin{proof}
First we check the properties (a) and (b) from Lemma \ref{crit} for both lists in the case where $n_1+n_2=1$.

We let $i\in\{1,2\}$ and $j:=3-i$. Since $d_{121}\neq 0 \neq d_{212}$ by hypothesis, we find elements $a\in C_i$, $b\in C_j$ and $g\in G$ in such a way that $a^g\cdot b=a$. Then  \[1=a^{-1}\cdot a^g\cdot b=[a,g]\cdot b.\] Therefore the tuple $((a),(g),(b))$ satisfies (a) and (b) for the list $[G, g,1\mid [m_1,1]]$ if $i=2$, and for the list $[G,g,1 \mid [m_2,1]]$ if $i=1$.

If we look at the list
$[G,\bar{g},0 \mid [m_1,2], [m_2,1]]$ and the case $i=1$, then the tuple $((),(),(a^{-1},a^g,b))$ satisfies (a) and (b). 

Finally, for the list $[G,\bar{g},0 \mid [m_1,1],[m_2,2]]$ and for the case $i=2$, we see that the tuple 
$((),(),(b^{-1},(a^g)^{-1},a))$ satisfies (a) and (b) of Lemma \ref{crit}.\\

Next we look at the case where $n_1+n_2=2$.
Since $d_{121}\neq 0$, we find  $a,c_1\in C_1$ and $c_2\in C_2$ such that $c_1 \cdot c_2=a$. Also, we have $d_{212}\neq 0$ and hence there exist $x\in C_2$ and $g\in G$ such that $x^g \cdot a=x$ (and $x^{-1} \cdot x^g=a^{-1}$). We calculate that
\[1=a^{-1}c_1c_2=x^{-1}x^g c_1c_2=[x,g]c_1c_2,\] showing that the tuple $((x),(g),(c_1,c_2))$ satisfies (a) and (b) from Lemma \ref{crit} for the list $[G,g,1 \mid [m_1,1],[ m_2,1]]$.

Now suppose that every maximal subgroup $M$ of $G$ intersects $C_1$ or $C_2$ trivially and take the tuples from above for the corresponding lists as stated. 
Set $U:=\langle a,b,g\rangle$ in the case where $n_1+n_2=1$
(and $U:=\langle a^{-1}, a^b,c\rangle=\langle c^{-1},(a^b)^{-1},a\rangle$ for the list $\bar{l}$), and set $U:=\langle x,g,c_1,c_2\rangle$ if $n_1+n_2=2$.
Then $U$ contains elements from $C_1$ and $C_2$, which forces $U=G$ and hence property (c) from Lemma \ref{crit}.
We argue similarly with the hypothesis about the orders of maximal subgroups.
\end{proof}

Since we have several cases for the group $\PSL_2(q)$, we 
need several lemmas and begin with a hypothesis for this case:

\begin{hyp}\label{PSL}
Suppose that $f\in\N$, that $p$ is an odd prime, that $q:=p^f\ge 9$ and $G=\PSL_2(q)$. Let
	$\epsilon\in\{1,-1\}$ be such that $q\equiv \epsilon$ modulo $4$ and set $\alpha:=\frac{q-\eps}4$, $\beta:=\frac{q-\eps}2$ and $\gamma:=\frac{q+\eps}2$.
\end{hyp}

\begin{lemma}\label{PSL21}
  If Hypothesis \ref{PSL} holds, where $q\notin\{9,11,19\}$, and if $l:=[G,g,1 \mid [\alpha, 1]]$ is a Hurwitz datum, then $l$ has a witness.    
\end{lemma}

\begin{proof}
	Let $c_{1,1}\in G$ be such that $o(c_{1,1})=\alpha$ and take a subgroup $I\le G$ of order $\beta$ that contains $c_{1,1}$. Let $t$ denote the central involution in the dihedral group $N:=N_G(\langle c_{1,1}\rangle)=N_G(I)$ and let $x,y \in G$ be involutions such that $N=\langle x,y\rangle$. We note that $N=C_G(t)$ has order $4 \cdot\alpha$. 
Some calculations show that $o(xy)=2 \cdot \alpha$ and that $[x,y]=(xy)^2$ has order $\alpha$. Hence we may choose $x,y$ such that $[x,y]={c_{1,1}}^{-1}$.
In particular $x,y \notin Z(N)$ and we can choose  $z\in C_G(x)\setminus N$. Then, if we set $a_1:=x$ and $b_1:=zy$, we see that ${c_{1,1}}^{-1}=[x,y]=[x,zy]=[a_1,b_1]$ and that $((a_1),(b_1),(c_{1,1}))$ satisfies Lemma \ref{crit}(a) and (b) for the list $l$. We turn to (c) and let $T:=\langle a_1,b_1,c_{1,1}\rangle$. 

Assume for a contradiction that $T\neq G$.
Using the subgroup structure of $G$ (see for example Dickson’s Theorem, as in Theorem
6.5.1 of [5]), we narrow down the possibilities for  
$T$ and $q$:

$T$ contains a $p'$-element and several involutions, therefore it is neither an elementary abelian $p$-group nor cyclic. It is also not contained in the normaliser of a Sylow $p$-subgroup, because $T'$ is not a $p$-group. 
If $T$ is a dihedral group, then $T=N$ and therefore $y \in T$ and finally $z \in T=N$, contrary to our choice of $z$. 

Next, $T$ could be isomorphic to a subgroup of $\Alt_4$, $\Sym_4$ or $\Alt_5$, with restrictions on $q$. 

Assume that $T$ is isomorphic to a subgroup of $\Alt_4$. Then $c_{1,1}\in T'$ forces $\alpha=2$ and hence $q\in\{7,9\}$, contrary to our hypothesis. 

Now assume that $T$ is isomorphic to a subgroup of $\Sym_4$. In particular $q^2 \equiv 1$ modulo $16$ for this to be possible, and we recall that $q>9$ and $q\notin\{11,19\}$. Then the element orders in $\Sym_4$ force $\alpha\in\{3,4\}$. If $\alpha=3$, then $q=13$, which does not work with the above congruence. 
Hence $\alpha=4$. But as before $c_{1,1}\in T'$, which is impossible now.  

Finally, assume that $T$ is isomorphic to a subgroup of $\Alt_5$. Then $q^2\equiv 1$ modulo $5$, and the possible element orders in $\Alt_5$ force $\alpha\in\{3,5\}$. If $\alpha=3$, then $q= 13$, which is impossible by the above congruence. Hence $\alpha=5$ and $q=19$, which is excluded by our hypothesis. 

Finally, assume that $T$ is contained in a subfield subgroup. First we let $m$ be a proper divisor of $f$ such that $T$ is isomorphic to a subgroup of $\PSL_2(p^m)$. Let $k,m\in\N$ be such that $f=k \cdot m$. We recall that $\alpha$ is coprime to $p$ and therefore $\alpha$ divides $\frac{p^m-1}{2}$ or $\frac{p^m+1}{2}$.  
This leads to $p=3$, $m=1$ and $f=2$, whence $q=9$, which is a contradiction.

If $k,m\in\N$ are such that $f=2km$ and $T$ is isomorphic to a subgroup of $\PGL_2(p^m)$, then we argue similarly and we deduce that $\alpha\le p^m+1$ and $q=25$. Now $T$ is isomorphic to a subgroup of $\PGL_2(5)\cong\Sym_5$ and $\alpha=6$. By construction  $T$ contains a dihedral group of order 24 and hence an element of order 12. But this is impossible. 

We conclude that $T=G$ and that $l$ has a witness, by Lemma \ref{crit}.
\end{proof}

\begin{lemma}\label{PSL22}

Suppose that Hypothesis \ref{PSL} holds and that
$q=9$. If $(g_0,n)\in\{(1,0),(1,1),(0,2)\}$ and $l:=[G,g,g_0 \mid [\alpha, 1],[\gamma, n]]$ is a Hurwitz datum, then $l$ has a witness. 
\end{lemma}

\begin{proof}
Since $q=9$, we know that $\epsilon=1$, $\alpha=2$ and $\gamma = 5$. 
Let $C_1$, $C_2$ denote a conjugacy class of involutions and of elements of order 5 in $G$, respectively. We calculate $d_{121}=8$ and $d_{212}=10$, and using Lemma \ref{max} we find elements $a_1,x_2,x_3,y_3\in C_2$, $a_2,b_1,b_2\in G$ and $c_1,c_2,c_3\in C_1$ such that the following hold: 

$((a_1),(b_1),(c_1))$ satisfies Lemma \ref{crit}(a) and (b) for the list $l$ if $(g_0,n)=(1,0)$,
$((a_2),(b_2),(c_2,x_2))$ satisfies the same criteria if $(g_0,n)=(1,1)$ and \\$((),(),(c_3,x_3,y_3))$ satisfies them if $(g_0,n)=(0,2)$. 

In each case we find a generating tuple $(x,y,z)$ for $G$ that satisfies all criteria from Lemma \ref{crit} (using the subgroup structure of $\PSL_2(9)$), which gives that the list has a witness in each case:

If $(g_0,n)=(1,0)$, then we choose $x:=c_1$, $y:=a_1^{-1}$ and $z:={a_1}^{b_1}$.
If $(g_0,n)=(1,1)$, then we choose 
$x:=c_2$, $y:=x_2$ and $z:=[a_2,b_2]$, and finally, if
$(g_0,n)=(0,2)$, then we take $x:=c_3$, $y:=x_3$ and $z:=y_3$.      
\end{proof}

\begin{lemma}\label{PSL23}

Suppose that Hypothesis \ref{PSL} holds, that $q=11$ and that $l:=[G,g, 1 \mid [\alpha, 1]]$ is a Hurwitz datum. Then $l$ has a witness. 
\end{lemma}

\begin{proof}
Let $C_1$ denote a conjugacy class of $G$ of elements of order $\alpha=3$ and let $C_2$ denote the unique conjugacy class of $G$ of elements of order 6. In particular $C_2$ is closed under inverses. We let $c \in G$ be such that $c^{-1}\in C_1$, we calculate $d_{221}=13$ and then we find $x,y\in C_2$ such that $xy=c^{-1}$. Since $y$ and $x^{-1}$ are conjugate in $G$ and $d_{221}>1$, we also find $a,b\in G$ such that $a:=x^{-1}$ and $y=a^b$.
Now \[1=xyc=a^{-1}a^b c=[a,b]c,\] which shows that Lemma \ref{crit}(a) and (b) hold for the list $l$. 

Assume for a contradiction that $U:=\langle a,b,c\rangle\neq G$. 
Since $U$ contains an element of order 6, it must be contained in a maximal subgroup of structure $D_{12}$.
This forces $c\in\langle a\rangle$, contrary to the choice of $c$. Now $U=G$ and Lemma \ref{crit}(c) is satisfied, which gives that $l$ has a witness. 
\end{proof}

\begin{lemma}\label{PSL24}
Suppose that Hypothesis \ref{PSL} holds, that $q=19$ and that $l:=[G,g, 1 \mid [\alpha,1]]$ is a Hurwitz datum. Then $l$ has a witness.
    
\end{lemma}

\begin{proof}
Let $C_1, C_2$ denote conjugacy classes in $G$ of elements of order $\alpha=5$ and of order 9, respectively. We calculate $d_{212}=36$ and 
	$d_{121}=40$, which shows that Lemma \ref{max} is applicable. Checking the list of maximal subgroup of $G$ gives that there is no maximal subgroup with order divisible by $45$, the lowest common multiple of $5$ and $9$. Then the lemma tells us that $l$ has a witness.     
\end{proof}



\begin{lemma}\label{PSL26}
Suppose that Hypothesis \ref{PSL} holds and that $q>9$. If $(g_0,n)\in\{(1,1),(0,2)\}$ and if
$l:=[G,g,g_0 \mid [\alpha, 1],[\gamma,n]]$ is a Hurwitz datum, then it has a witness. 
\end{lemma}

\begin{proof}
\textbf{Case 1:} $(g_0,n)=(0,2)$.\\
Let $C_1$, $C_2$ denote conjugacy classes in $G$ of elements of order $\alpha$ and $\gamma$, respectively. Calculations in \cite{chev} give that $d_{212}=2q+2\epsilon$, where $\epsilon$ is as in our main hypothesis. Hence let $x\in C_1, y,z\in C_2$ be such that $yx=z$. Then $1=yxz^{-1}=xz^{-1}y$, whence the tuple $((),(),(x,z^{-1},y))$ satisfies Lemma \ref{crit}(a) and (b). Let $U:=\langle x,y,z\rangle$ and assume for a contradiction that $U\neq G$. Since $q>9$, we find elements of order $\alpha\ge 3$ and $\gamma\ge 5$ in $U$, and moreover $|U|$ is divisible by the lowest common multiple of these orders, i.e. $\frac{q^2-1}8$.
Then the subgroup structure of $G$ (see for example Theorem 6.5.1 of \cite{GLS3}) implies that 
$q=11$ and that $U\cong\Alt_5$, because $\frac{q^2-1}8=15$.
We let $C_{x}:=C_1\cap U$ and $C_y:=C_2\cap U$ and we see that these are conjugacy classes in $U$. There are $\frac{\abs G}{\abs U}=11$ conjugates of the subgroup $U$ in $G$ and moreover \[\abs{C_2}=\frac{\abs G}{\abs{ C_G(y)}}=2^2\cdot 3\cdot 11\qquad\text{and}\qquad\abs{C_{y}}=\frac{\abs U}{\abs{C_U(y)}}=2^2\cdot 3.\] 

Since $\abs{C_2}=\abs{C_{y}}\cdot\frac{\abs G}{\abs U}$, it follows that $U$ is the unique $G$-conjugate of $U$ containing $y$. We calculate $d_{y x y}$ with respect to conjugacy in $U$ and we obtain the value $5<20=d_{212}$. Then we choose $(a,b)\in C_2\times C_1$ such that $ab=z$ and $\langle a,b,z\rangle$ is not contained in any $G$-conjugate of $U$. Since $G$ has another conjugacy class of subgroups isomorphic to $U$, we can choose from the  $d_{212}=20$ pairs in $M:=\{(c,d)\in C_2\times C_1 \mid cd=z\}$ and we find exactly ten pairs that generate $G$.
For such a pair (c,d), the tuple 
tuple $((),(),(c,z^{-1},d))$ satisfies all criteria from Lemma \ref{crit}. Then $l$ has a witness in this case. 

\textbf{Case 2:} $(g_0,n)=(1,1)$.\\
Let $\bar{g}\in\N$ be such that $\bar{g}>1$ and $\bar {l}:=[G,\bar{g},0\mid [\alpha,1], [\gamma,2]]$ is a Hurwitz datum. By Case 1 we know that $\bar{l}$ has  a witness, so we let $((),(),(x,y,z))$ be such that the criteria from Lemma \ref{crit} are satisfied. 
Then $o(x)=\alpha$, $o(y)=\gamma=o(z)$, $1=xyz$ and $\langle x,y,z\rangle=G$. In particular $xy=z^{-1}$ and $1=zxy$.
Let $C_1,C_2$ denote conjugacy classes of $G$ of elements of order $\beta$ and $\gamma$, respectively, and without loss $z \in C_2$. 
Calculations in \cite{chev} show that $d_{121}=2q-2\eps$, and then we let $a,b\in C_1$ be such that $bz=a$. Moreover, let $g\in G$ be such that $b=a^g$. Then $1=a^{-1}a^g z=[a,g]z$ and therefore $1=zxy=[g,a]xy$. Also, $z\in\langle a,g\rangle$ and this implies that $G=\langle x,y,z\rangle= \langle x,y,a,g\rangle$. In particular $((g),(a),(x,y))$ satisfies the criteria from Lemma \ref{crit} and therefore $l$ has a witness in this case.
\end{proof}

\begin{lemma}\label{sgp1}
	Suppose that $p$ is a prime, $f\in\N$, let $q:=p^f\ge 4$, $d:=(q-1,2)$ and $\alpha:=\frac{q+1}d$. If $p=2$, then suppose that $f$ is even. Let $H:=\PSL_2(q)$. Then there are elements $x,y,z\in H$ that generate $H$ and such that $o(x)=2$, $o(z)=\alpha$ and $[x,y]=z$.
\end{lemma}

\begin{proof}
We use the subgroup structure of $H$, see for example Theorem 6.5.1 in \cite{GLS3}. 
Let $z\in H$ be an element of order $\alpha$. Then $N_H(\langle z\rangle)$ is dihedral of order $2\alpha$. If we let $C_1$ denote a conjugacy class of involutions and $C_2:=z^H$, then we can calculate (for example with \cite{chev}) that the set $M:=\{(a,b)\in C_1\times H \mid [a,b]=z^{-1}\}$ has exactly $\alpha\cdot\abs{C_H(t)}$ elements.
Now we find $(t,b)\in M$ such that $\langle t,b,z\rangle=H$:
If $q\in\{4,5,7,9\}$, then we can check this using \texttt{GAP} (\cite{GAP4}). 
If $q\notin\{4,5,7,9\}$, then $\alpha>5$ and $U:=N_H(\langle z\rangle)$ has order at least $12$ and is a maximal subgroup of $H$. Since $[t,b]=1\neq z$, we see that $t$ is not central in $U$, but that it inverts $z$. Now $\abs{C_H(t)\cap U}\in\{2,4\}$.
	Moreover $\abs{C_H(t)}\ge 8$ because $q \ge 11$, and we choose $g\in C_H(t)\setminus U$. Then  $(t,g\cdot b)\in M$ and we calculate that $V:=\langle t,g\cdot b,z\rangle$ is in fact all of $H$. 
\end{proof}

\begin{lemma}\label{sgp2}
Suppose that $p$ is prime, $f\in\N$ and let $q:=p^f\ge 4$, $d:=(q-1,2)$ and $\alpha:=\frac{q+1}d$. If $p=2$, then suppose that $f$ is even. Suppose that $G\cong\PSL_2(q).C_2$. Then there are elements $x,y,z\in G$ such that $o(x)=2$, $o(z)=\alpha$, $[x,y]=z$ and $\langle x,y,z\rangle=G$ and such that the following is true:
	\[\abs{C_G(x)}=\begin{cases} 2q,&\text{if }p=2,\\ 2(q-1),&\text{if } q\equiv 1 \mod 4,\\ 2(q+1),&\text{if }q\equiv -1\mod 4.\end{cases}\]
\end{lemma}

\begin{proof}
Let $H\le G$ be such that $H\cong\PSL_2(q)$ and, using Lemma \ref{sgp1} choose elements 
$x,y,z\in H$ such that $o(x)=2$, $o(z)=\alpha$, $[x,y]=z$ and $\langle x,y,z\rangle=H$. We prove that $\abs{C_G(x):C_H(x)}=2$ and then deduce our statement about $\abs{C_G(x)}$.\par\medskip

First suppose that $p\neq 2$. Then $\frac{q+1}2=\alpha$ or $\frac{q-1}2$ is even and $C_H(x)$ is dihedral of order $q+1$ or $q-1$. If $u\in C_H(x)$ is an element of order $\alpha$, then $N_H(\langle u\rangle)=C_H(x)$.Let $U:=\langle g\rangle$ and $\Delta:=U^G$. Then all elements of $\Delta$ are subgroups of $H$, because $H\unlhd G$, and then $H$ acts transitively on $\Delta$ by conjugation. If $V\in\Delta$, then its point stabiliser in this action is  $N_G(V)$ and Frattini gives $G=H\cdot N_G(U)$.
Moreover
	\[N_G(U)/N_H(U)\cong (H\cdot N_G(U))/H=G/H\]

 and therefore  $\abs{C_G(x):U}=\abs{N_G(U):U}=4$.
 Now $|C_G(x):C_H(x)|=2$.

 Next suppose that
$p=2$. Then $S:=C_H(x)$ is elementary abelian of order $q$ and it is a Sylow 2-subgroup of $H$. Moreover $H$ acts transitively on the set of subgroups of order 2 of $S$. Since $H\unlhd G$, we know that all conjugates of $x$ in $G$ are contained in $H$ and Frattini gives that $G=H\cdot N_G(\langle x\rangle)=H\cdot S$. Again we see that $|C_G(x):C_H(x)|=2$.

For all $p$, we use that $\abs{C_G(x):C_H(x)}=2$ and we choose $g\in C_G(x)\setminus H$ and set $w:=gy$. 
Then $1=[x,y]z=[x,w]z$ and we assume for a contradiction that $U:=\langle x,w,z\rangle\neq G$.
Since $w \notin H$, this assumption forces 
$U\neq H$.
Moreover $U\neq \langle x,z\rangle$, because $\langle x,z\rangle\ug H$. Let $a\in U\setminus H$. 
Since $G/H=\{H,Ha\}$, we see that $a\cdot w=agy\in U\cap H$, and we may suppose that $agy\notin \langle x,z\rangle$.
Then we choose $a$ such that $H=\langle x,z,agy\rangle\ug U$ and it follows that $y\in \langle x,z,agy\rangle=H$ and $y\in U$. Then $g\in U$, contrary to our assumption.
\end{proof}

One more remark is due before we turn to the series
$\PSp_4(q)$. 

\begin{rem}\label{mat}
Let $q$ be a power of $2$, $q \ge 4$. We investigate a particular involution centraliser in $\PSp_4(q)$, referring to \cite{AschbacherSeitz76} and its notation on pages 8 -- 10 and 14 -- 24.

Let $f\in\N$, $f\ge 2$, $q:=2^f$, $n=4$ and let $V$ denote an $n$-dimensional vector space over
	$\GF(q)$. Moreover, let $G:= \Sp(V)\cong\PSp(V)$ and let $t\in G$ be an involution. We describe $C_G(t)$ in the case where $t$ is as in  (7.11) in \cite{AschbacherSeitz76}, using matrices. We denote the unity matrix of size $m \times m$ ($m \in \N$) by $I_m$, we use $*$ for the multiplication of matrices and we omit zero entries. 

Following the hypotheses in (4.2) and (7.6) of \cite{AschbacherSeitz76}, we have that $1\le l\le \frac n2=2$ is even, hence $l=2$. We set $C_l:=C_{\SL(V)}(t)$ and $C:=C_l\cap\Sp(V)=C_G(t)$. Then (4.2) in \cite[p. 9]{AschbacherSeitz76} gives that we can choose $t$  to be \[t=j_l=\begin{pmatrix}
	I_l\\&I_{n-2l}\\I_l&&I_l
	\end{pmatrix}.\] 
 Let $g \in C_l$. Then (4.3) in \cite{AschbacherSeitz76}
 yields that
	\[g=\begin{pmatrix}
	X(g)\\
	P(g)&Y(g)\\
	Q(g)&R(g)&X(g)
	\end{pmatrix},\] 

 where 
$X(g)$ and $Q(g)$ have size $l\times l$, $Y(g)$ has size $(n-2l)\times(n-2l)$, $P(g)$ has size $(n-2l)\times l$ and $R(g)$ has size $l\times(n-2l)$. Moreover $\det(X)^2\det(Y)=1$. Since $l=2$, the matrices $Y(g),P(g)$ and $R(g)$ do not occur for all $g\in C$.

Let $g\in C$. Then we use (7.11) in \cite{AschbacherSeitz76}, and it follows (with their notation) that \[F=\begin{pmatrix}
	0&1\\1&1
	\end{pmatrix}, X(g)=\begin{pmatrix}
	1&0\\x&1
	\end{pmatrix}, Q(g)=\begin{pmatrix}
	r&z\\y&s
	\end{pmatrix} \text{ and } g=\begin{pmatrix}
	X(g)\\Q(g)&X(g)
	\end{pmatrix},\]
	where $x,y,z,r,s\in \GF(q)$.

In addition, we see that \[\begin{pmatrix}
	0&0\\0&0
	\end{pmatrix}=X(g)*F*Q(g)^T+Q(G)*F*X(g)^T,\] using (7.11), where $M^T$ denotes the transpose of the matrix $M$. Multiplication gives \[M_1:=X(g)*F*Q(g)^T=\begin{pmatrix}
	z&s\\r+z(x+1)&y+s(x+1)
	\end{pmatrix}\]
 and 
\[M_2:=Q(g)*F*X(g)^T=\begin{pmatrix}
	z&z(x+1)+r\\s&s(x+1)+y
	\end{pmatrix}.\] We recall that $\GF(q)$ has characteristic $2$, and then (7.11) in \cite{AschbacherSeitz76} implies that \[\begin{pmatrix}
	0&0\\0&0
	\end{pmatrix}=M_1+M_2=\begin{pmatrix}
	0&s+z(x+1)+r\\
	s+z(x+1)+r&0
	\end{pmatrix}.\]
	Now if we set $s=-r-z(x+1)=r+z(x+1)$, the we can choose $y,r,z,x\in F$ as we like. 
 
Therefore \[C=\{\begin{pmatrix}
	1&0&0&0\\
	x&1&0&0\\
	r&z&1&0\\
	y&r+z(x+1)&x&1
	\end{pmatrix} \mid x,y,z,r\in \GF(q)\}.\]
	We apply Satz II.7.1 in \cite{Huppert-I} to see that $C$ is a subgroup of a Sylow 2-subgroup~$S$ of $\SL(V)$. Here $\abs S=q^6$ and the matrices of $S$ have arbitrary entries below the main diagonal. We compare this with the entries of matrices from $C$ and the relationships between them, and then it follows that $|C|=\frac{\abs S}{q^2}=q^4$. We compare this to $\abs G$ (see for example (3.22) in \cite[p. 60]{Wilson}) and see that $C \in \Syl_2(G)$.
\end{rem}

 \begin{lemma}\label{Sp4}
	Suppose that $q\ge 3$ is a prime power, that $d:=(q^2-1,2)$, $\alpha:=\frac{q^2+1}{d}$ and $G=\PSp_4(q)$. If $l:=[G, g,1 \mid [\alpha, 1]]$ is a Hurwitz datum, then it has a witness. 
\end{lemma}

 \begin{proof}
Theorems 3.7 and 3.8 in \cite[pp. 92]{Wilson} tell us that $G$ has a maximal subgroup $M$ of structure $\PSp_2(q^2).C_2\cong\PSL_2(q^2).C_2$. Then Lemma \ref{sgp2} gives elements $x,y,z\in M$ that generate $M$ and such that $o(x)=2$, $o(z)=\alpha$ and $[x,y]z=1$.

If $q$ is even, then $\abs{C_M(x)}=2q^2$ by Lemma \ref{sgp2}, and if $q$ is odd, then $q^2\equiv 1$ modulo  $4$ and $\abs{C_M(x)}=2(q^2-1)$ by the same lemma. 

If $q$ is odd, then $C_G(x)$ has a section isomorphic to the central product $\SL_2(q)\circ \SL_2(q)$ or to $\PSL_2(q)$, by Table 4.5.1 in \cite{GLS3}. 

If $q$ is even, then $C_G(x)$ has a section isomorphic to $\PSL_2(q)$ by (7.9) and (7.10) in \cite{AschbacherSeitz76} or $C_G(x)$ is a 2-group by (7.11) in \cite{AschbacherSeitz76}. In the case of (7.11), it is also true that $C_G(x)$ is a Sylow 2-subgroup of $G$, with order $q^4$. In the first case,  (7.9) or (7.10), it follows that $\abs{C_G(x)}\ge q^2(q^2+1)> q^4$, because $C_M(x)\le C_G(x)$ and $\abs{C_M(x)}=2q^2$.
In both cases we conclude that 
$\abs{C_M(x)}<\abs{C_G(x)}$. 

Now we look at Lemma \ref{crit} (a) and (b):

We recall that $x,y,z \in M$ generate $M$ with additional properties, and if we take $g\in C_G(x)\setminus M$, then  $[x,gy]z=[x,y]z=1$. Therefore the tuple $((x),(gy),(z))$ satisfies Lemma \ref{crit} (a) and (b)for the list $l$.
	
 Let $T:=\langle x,g\cdot y,z\rangle$. If $T=G$, 
 then Lemma \ref{crit} (c) holds and $l$ has a witness. 
Otherwise $T$ is contained in a maximal subgroup $U$ of $G$, and $\langle x,z\rangle\le M$. Since $y\in M$ and $g\notin M$, we deduce that $T\nleq M$.

Suppose that $q$ is odd. Then 
$\alpha=\frac{q^2+1}2$ and $2\alpha$ divides $\abs T$ as well as $\abs U$. If $q\ge 5$, then $\alpha\ge 13$ and we apply Theorem 3.8 of \cite{Wilson} 
together with formulae for the orders of classical groups:
Then $U\cong M\cong\PSL_2(q^2).C_2$ or $q=3$, $\alpha=5$ and $U\cong {C_2}^4.\Omega_4^-(2)\cong {C_2}^4.\Alt_5$, 
checking the maximal subgroups from 
Theorem 3.8 in \cite[p. 93]{Wilson} for $\alpha$ dividing their order. 
Now suppose that $q=3$. We let $C_1$ denote a conjugacy class of $G$ of elements of order $\alpha=5$ and we let $C_2$ denote a conjugacy class of elements of order $9$. Then we calculate $d_{121}=585$ and $d_{212}=567$, and we see that proper subgroup of $G$ cannot intersect both $C_1$ and $C_2$ trivially. Then Lemma \ref{max} gives that $l$ has a witness. 

Finally, suppose that $q>3$. Then $U\cong M$ (see above) and Table 3.5.C of \cite[p. 72]{KleidmanLiebeck} (Block $C_3$, where $j=10$, column V, explained on pp. 61) shows that $M^G$ is the unique class of maximal subgroups of $G$ that are isomorphic to $M$. Thus let $h\in G$ be such that $U^h=M$. Then $\langle z\rangle$ is a Hall $\pi(o(z))$-subgroup of $G$ and $N_G(\langle z\rangle)\le M$ (see for example Theorems 3.7 and 3.8 in \cite{Wilson}). 

The results III.5.8 and III.5.6 in \cite{Huppert-I} give that we may choose $h\in G$ in such a way that $\langle z\rangle^h=\langle z\rangle$. Consequently $h\in N_G(\langle z\rangle)\le M$ and then $U=M$, contrary to our assumption.

Now suppose that $q$ is even. 
Then $\alpha=q^2+1$ and $2\alpha$ divides both $\abs T$ and $\abs U$. Since $q\ge 4$, we know that $\alpha\ge 17$. Then we apply Theorem 3.7 in \cite{Wilson} and order formulae for exceptional groups, and we deduce that $U \cong M\cong\SL_2(q^2).C_2$ or $U \cong \GO_4^-(q)\cong \SL_2(q^2).C_2$ or $f$ is odd and $U\cong \Sz(q)$.
 
With information about the structure of Suzuki groups (for example XI.3 in \cite{Huppert-III}) we can exclude the possibility $U\cong\Sz(q)$. Then $U\cong\SL_2(q^2).C_2$ and we use Table 3.5.C in \cite[p. 72]{KleidmanLiebeck}, which gives that $U$ is conjugate to $M$ or it lies in the other class of maximal subgroups isomorphic to $M$ (see block $C_3$, where $j=10$ and block $C_8$ from Table~3.5.C in \cite[p. 72]{KleidmanLiebeck}). 
The first case gives a contradiction, as above. 
Now $U\cdot M\subseteq G$ and \[\abs U\cdot \abs M=\abs U^2=4q^4(q^2-1)^2(q^2+1)^2>q^4(q^2-1)^2(q^2+1)=\abs G.\]

 Therefore $\abs{U\cap M}=\frac{\abs U\cdot \abs M}{\abs{U\cdot M}}\ge 4(q^2+1)=4\alpha$ and $\abs{U\cap M}$ is a multiple of $4(q^2+1)$. Hauptsatz II.8.27 from \cite{Huppert-I}, together with the fact that $q\ge 4$, gives that $\abs{U\cap M}=4\alpha$ and $U\cap M=N_G(\langle z\rangle)$. 
 
Moreover $\abs{C_U(x)}=2q^2=\abs{C_M(x)}$ and $\abs{C_{U\cap M}(x)}\le 4$, since $o(x)=2$ and $\alpha$ is odd. We had already seen earlier that $\abs{C_G(x)}\ge q^4$, and then $\abs{C_G(x)\setminus(U\cup M)}\ge 4$. Choosing $g\in C_G(x)\setminus M$ as above, such that $g y\notin U\cup M$, it follows that $T=G$.

In all cases that we looked at, Lemma \ref{crit} gives that $l$ has a witness. 
 \end{proof}

 \begin{lemma}\label{POmega}
     Suppose that $q\ge 2$ is a prime power, set $d:=(q^4-1,2)$ and $\alpha:=\frac{q^4+1}{d}$ and suppose that $G=\POm_8^-(q)$. If $l:=[G, g,1 \mid [\alpha, 1]]$ is a Hurwitz datum, then $l$ has a witness.
 \end{lemma}

\begin{proof}
	By Theorem 3.11 in \cite[pp. 94]{Wilson}, $G$ has a  maximal subgroup $M$ of structure $\POm_4^-(q^2).C_2\cong\PSL_2(q^4).C_2$. Then we use Lemma \ref{sgp2} and we find $x,y,z\in M$ that generate $M$ and such that $o(x)=2$, $o(z)=\alpha$ and $[x,y]z=1$.
	If $q$ is even, then the same lemma tells us that $\abs{C_M(x)}=2q^4$. If $q$ is odd, then $q^4\equiv 1$ modulo $4$ and $\abs{C_M(x)}=2(q^4-1)$. Moreover 
 $C_G(x)$ has a section isomorphic to $\Omega_6^-(q)$, to $\Omega_6^+(q)$ or to $\Omega_4^-(q)\times \Omega_4^+(q)$ in this case, as we can see in Table 4.5.1 in \cite{GLS3}. All these sections have order divisible by  $q$. 
 
 Returning to the case where $q$ is even, we see that $C_G(x)$ has a section isomorphic to \mbox{$\Sp_2(q)\times \POm_4^-(q)$}, to $\Sp_4(q)$ or to $\Sp_2(q)\cong \SL_2(q)$, by (8.6), (8.8), (8.10) and (7.6) in \cite{AschbacherSeitz76}. In particular these sections have order divisible by $q^2-1$. 
	Since $\abs{C_G(x)}>\abs{C_M(x)}$, we find an element $g\in C_G(x)\setminus M$ and then $[x,gy]z=[x,y]z=1$. The tuple $((x),(gy),(z))$ satisfies Lemma \ref{crit}(a) and (b) now, and we only need to check (c). 
 Set $T:=\langle x,g y,z\rangle$ and assume for a contradiction that $T$ is contained in a maximal subgroup $U$ of $G$. We recall that $\langle x,z\rangle\le M$. Since $y\in M$ and $g\notin M$, it follows that $T\nleq M$.
 
Now we use that $G$ is a section of $\GO_8^-(q)$ and that we have a list of maximal subgroups in Theorem 3.11 in \cite[pp. 94]{Wilson}. Let $\hat M$ denote a maximal subgroup of $\GO_8^-(q)$ such that $U$ is a section of $\hat M$. Since $T\le U$, we see that $\abs T=\alpha=d^{-1}(q^4+1)$ divides $|\hat M|$. Also, we can exclude the cases (v) to (vii), (x) and (xi) from Theorem 3.11. Moreover (i) to (iv) and (viii) cannot occur because, in these cases, $\alpha$ would not divide $|\hat M|$.
The final remaining case is $\hat M\cong \GO_4^-(q^2).C_2$ and $U\cong \POm_4^-(q^2).C_2\cong\PSL_2(q^4).C_2\cong M$.

 Using Table 3.5.F in \cite[p. 74]{KleidmanLiebeck} (block $C_3$, where $j=16$, column V), we see that $M^G$ is the unique class of subgroups of $G$ isomorphic to $M$. For an explanation for column V we refer to 3.2 in \cite[pp.~61]{KleidmanLiebeck}. Let $h\in G$ be such that $U=M^h$. There is a unique $G$-conjugacy class of cyclic subgroups of order $\alpha$ in $G$ and $N_G(\langle z\rangle)\le M$. Moreover $T$ is a Hall subgroup of $G$, because $\alpha=\frac{q^4+1}d$ is coprime to $q$, to $q^4-1$, to $q^6-1$ and hence to ~$\frac{\abs G}{\alpha}$. The results III.5.8 and III.5.6 in \cite{Huppert-I} imply that we may choose $h\in G$ in such a way that $\langle z\rangle^h=\langle z\rangle$, and then $h\in N_G(\langle z\rangle)\le M$. But this forces $U=M$, which is a contradiction.
    Hence $l$ has a witness by Lemma \ref{crit}.
\end{proof}

\begin{lemma}\label{suz}
Suppose that $f\in\N$, let $s:=2^f$, $q:=2s^2\ge 8$ and $G=\Sz(q)$, moreover let $\alpha:=q-2s+1$, $\beta:=q-1$ and
	$\gamma:=q+2s+1$. Then $\alpha$, $\beta$ and $\gamma$ are pair-wise coprime. 
	
  Let $n_1,n_3\in\{0,1\}$ and $n_2\in\{0,1,2\}$ be such that
	$1\le n_1+n_2+n_3\le 4$.
If $n_1+n_2+n_3\ge 3$, then set $g_0=0$. Otherwise $g_0=1$.

 Then, if $l := [G,g,g_0 \mid 
 [\alpha,n_1], [\beta,n_2], [\gamma,n_3]]$ is a Hurwitz datum, it has a witness. 
\end{lemma}

\begin{proof}

Let $C_1$, $C_2$ and $C_3$ denote conjugacy classes of elements of order ~$\alpha$, $\beta$ and $\gamma$, respectively. We note that $\alpha\cdot\gamma=q^2+1$. Since $\beta$ and $q^2+1$ are coprime, it follows that $(\alpha,\beta)=1=(\beta,\gamma)$. Moreover $(\alpha,\gamma)=1$ and $\alpha\ge 5$, $\beta\ge 7$ and $\gamma\ge 13$, and all these numbers are odd. Then we deduce that $\alpha$, $\beta$ and $\gamma$ are all coprime to $q$. 

Next we turn to \cite{chev} and we see that $C_1$ belongs to the generic conjugacy class 7 in \cite{chev}, $C_2$ belongs to class 5 and $C_3$ belongs to class 6.

We prove a little statement (U):

Suppose that $i,j\in\{1,2,3\}$ are distinct and that $U\le G$ intersects $C_i$ and $C_j$ non-trivially. Then we prove that $U=G$:

Assume otherwise and let $M$ be a maximal subgroup of $G$ containing $U$. 
By our hypothesis in (U), both $\abs U$ and $\abs M$ are divisible by 
	$\alpha\cdot \beta=(q-1)(q-2s+1)$, by $\alpha\cdot\gamma=q^2+1$ or by $\beta\cdot\gamma=(q-1)(q+2s+1)$.
	Then Theorem 4.1 in \cite[p. 117]{Wilson} (if $q>8$) and \cite[p. 28]{ATLAS} (if $q=8$) immediately give a contradiction. 

 Suppose that $n_1=1=n_3$ and $n_2\in\{1,2\}$.
 Using the generic character table of $G$ in \cite{chev}, we calculate: \[d_{123}=q^3+q^2(2s+1)+q+2s+1>0.\] 
 This gives elements $x\in C_1$, $y\in C_2$ and $z\in C_3$ such that $xy=z$. This gives a tuple that satisfies Lemma \ref{crit}(a) and (b). If $(n_1,n_2,n_3)=(1,1,1)$, then our claim (U) above implies that $l$ has a witness, with Lemma \ref{crit}. 
If $(n_1,n_2,n_3)=(1,2,1)$, then we can use Lemma 3.7 in \cite{SW} (because $\beta=q-1$ is odd) and it shows that 
$l$ has a witness again. 

For the remaining cases of $n_1,n_2$ and $n_3$, we 
recall that $s=2^f$ and we
calculate
	\begin{align*}
		d_{121}&=q^3+(1+2s)q^2+q+(1+2s), &d_{212}&=q^3+(1+2s)q^2-q-(1+2s),\\
		d_{232}&=q^3+(1-2s)q^2-q-(1-2s), &d_{323}&=q^3+(1-2s)q^2+q+(1-2s),\\
		d_{131}&=q^3-q^2+(4s-1)q+1\text{ and} &d_{313}&=q^3-q^2-(4s+1)q+1.\\
	\end{align*}
 Whenever $i,j\in\{1,2,3\}$ and $i\neq j$, then 
	$d_{iji}>0$ because $s>1$ and $q\ge 8$.
 
 Our claim (U) and Lemma \ref{max} give that $l$ has a witness in almost all cases, the only exception is 
$(n_1,n_2,n_3)=(0,2,0)$. 
Here we use the case $(n_1,n_2,n_3)=(0,1,0)$ and Lemma 
3.7 in \cite{SW} once more.    
\end{proof}

\begin{lemma}\label{d4}
	Suppose that $q\ge 2$ is a prime power, set $\alpha:=q^4-q^2+1$ and let $G=\, ^3D_4(q)$. If the list $l:=[G,g,1 \mid [\alpha, 1]]$ is a Hurwitz datum, then it has a witness. 
\end{lemma}
\begin{proof}
	Let $C_1$ denote a conjugacy class of elements of order $\alpha$. Theorem~4.3 (xi) in \cite[p. 144]{Wilson} gives us an element $x\in G$ of order $\beta:=q^2+q+1$ in such a way that $\abs{C_G(x)}$ is divisible by $\beta^2$ and that $C_G(x)$ lies in a maximal subgroup of $G$ of structure $(C_\beta\times C_\beta):\SL_2(3)$. Then we set $C_2:=x^G$.

First we prove that any subgroup $U$ of $G$ that contains elements of order $\alpha$ and $\beta$ must be $G$. 
If this were not the case, then $G$ would have a maximal subgroup $M$ of order divisible by 
$\alpha=q^4-q^2+1\ge 2^4-2^2+1=13$, and then the formulae for $\abs{\SL_2(q^3)}$, $\abs{\SL_2(q)}$, $\abs{G_2(q)}$, $\abs{\PGL_3^\pm(q)}$, $\abs{\PSL_2(q^3)}$, $\abs{\PSL_2(q)}$, $\abs{\SL_3^\pm(q)}$, $\abs{C_{\beta}}$, $\abs{C_{q^2-q+1}}$ and $\abs{\SL_2(3)}=8\cdot 3$ leave only one possibility: $M$ must have structure 
$C_\alpha:C_4$, by Theorem 4.3 in \cite[p. 144]{Wilson}. But $\beta=q^2+q+1\ge 2^2+2+1=7$ also divides $|M|$, which is impossible now. 

Using the generic character tables for $G$ in \cite{chev}, we see that $C_1$ belongs to the generic class 28, if $q$ is even, and to class 32, if $q$ is odd. $C_2$ belongs to the generic class 26, if $q$ is even, and to the class 30, if $q$ is odd. 
We calculate
	\begin{align*}
		d_{121}=q^{20}&-2q^{19}+q^{18}+2q^{17}-4q^{16}+2q^{15}+q^{14}-2q^{13}+q^{12}+4q^{11}-8q^{10}\\&
		+4q^9+q^8-2q^7+q^6+2q^5-4q^4+2q^3+q^2-2q+1 \qquad \text{and}\\
		d_{212}=q^{20}&-2q^{19}+q^{18}+2q^{17}-4q^{16}+2q^{15}+q^{14}-2q^{13}+q^{12}-8q^{11}+16q^{10}\\&
		-8q^9+q^8-2q^7+q^6+2q^5-4q^4+2q^3+q^2-2q+1.
	\end{align*}

 Then Lemma \ref{max} implies that $l$ has a witness. 
\end{proof}

\begin{lemma}\label{G2}
Suppose that $f\in \N$ and let $s:=3^f$, $q:=3s^2$, $\alpha:=\frac{q-1}2$ and $G=\, ^2G_2(q)$. If $l:=[G,g,1 \mid [\alpha, 1]]$ is a Hurwitz datum, then $l$ has a witness.    
\end{lemma}

\begin{proof}
We let $C_1$ denote a conjugacy class of elements of order $\alpha$. Using Theorem~4.2~(iv) in \cite[p. 138]{Wilson} we take $\beta:=q+3s+1$ and we let $C_2$ denote a conjugacy class of elements of order $\beta$.

Let $U \le G$ be such that $U$ contains elements of orders $\alpha$ and $\beta$. Then we prove that $U=G$:
Theorem 4.2 in \cite[p. 138]{Wilson} shows that every maximal subgroup of $G$ that has an element of order $\beta$ has structure $C_\beta:C_6$. This is because  $\beta$ divides neither $\abs{\PSL_2(q)}$ nor $\abs{^2G_2(q_0)}$ (if $q_0$ is a proper divisor of $q$).  Since $q\ge 27>6$ and $U$ also has an element of order $\alpha\ge 13$, it cannot be contained in a maximal subgroup of $G$. 

Now we do calculations with the generic character table of $G$ in \cite{chev}, where $C_1$ belongs to the generic class 11 and $C_2$ belongs to class 14. Then 
	\begin{align*}
		d_{121}&=q^5+(2-3s)q^4+(1-3s)q^3+q^2+(2-3s)q+(1-3s)&\text{and}\\
		d_{212}&=q^5+(2-3s)q^4+(1-3s)q^3-q^2-(2-3s)q-(1-3s).
	\end{align*}
	For all $s\ge 1$, we see that $d_{121}>0$, $d_{212}>0$ and hence Lemma \ref{max} gives that $l$ has a witness. 
    
\end{proof}

\section{Proof of the main result}

Suppose that $G$ is as described in the first hypothesis of Theorem \ref{main4fpH} and let $\Omega$ be a $G$-orbit of $X$ on which $G$ acts with fixity 4. Then Theorem 1.3 in \cite{BHMSW} gives the possibilities for $G$ (bearing in mind that the point stabilisers are cyclic),
which we have listed in Table \ref{GroupList} for convenience. The table also gives the order of the point stabiliser, respectively, and indicates whether or not there is a possibility for mixed fixity. If $H \le G$ is a point stabiliser in the action on $\Omega$, then every element in $H^\#$ fixes exactly four points in total, which implies that there is only one non-regular orbit of $X$ on which $G$ acts as it does on $\Omega$. 

In the cases without mixed fixity, there is only one non-regular orbit and hence a single possibility for the branching datum.  
This gives exactly the lists in Lines 6, 7, 9-11 and 13 of Table \ref{DataList}. 

Now we consider the cases with mixed fixity as given by Table \ref{GroupList}, one by one. For the orders of point stabilisers in the cases where the group acts with fixity 2 or 3, we refer to Lemmas 3.11 and 3.13 in \cite{MW} and Theorem 1.1 in \cite{MW3}. 

If $G=\Alt_7$ and there is only one non-regular orbit, then this gives Line 1a of Table \ref{DataList}. There is also exactly one possibility for a fixity 3 action, with point stabilisers of order $7$.
This means that $G$ can have a second non-regular orbit on $X$ on which it acts with fixity 3, but this can only happen once because of the global fixity 4 hypothesis. This gives the branching datum in Line 1b.

If $G=\PSL_2(7)$ and there is only one non-regular orbit, then this gives Line 2a of Table \ref{DataList}. There is also exactly one possibility for a fixity 2 action, with point stabilisers of order $3$. Under the global fixity 4 hypothesis it is possible for $G$ to act on up to two orbits with fixity 2 and with the same point stabilisers, and this gives the lines 2b and 2c. 
The group can also act with fixity 3 and point stabilisers of order 7, but then this can only happen on a single orbit. This gives rise to line 2d, with fixity 3 and 4 combined, and if we allow $G$ to act with fixity 2, 3 and 4 on different orbits, then this gives the branching data in Lines 2e and 2f.

In the case where $G=\PSL_2(8)$, we have a similar situation.
A unique non-regular orbit with fixity 4 action gives Line 3a of Table \ref{DataList}. A combination with a fixity 2 action, for which we have two possible point stabilisers, leads to up to five non-regular orbits and hence to Lines 3b--3i. 

Now suppose that $G=\PSU_4(3)$.
If there is a unique non-regular orbit with fixity 4 action, then this gives Line 4a of Table \ref{DataList}. In the mixed case with fixity 3 and 4, we obtain Line 4b.

If $q$ is a prime power, $q \ge 9$, $q$ is odd and $G=\PSL_2(q)$, then following Table \ref{GroupList} we
see two possibilities for a fixity 4 action, and if $G$ has only one non-regular orbit, then this gives Lines 5a and 5b
of Table \ref{DataList}.
In the mixed fixity cases, we first suppose that 
$q \equiv 1$ modulo $4$. Then there can be up to two orbits where $G$ acts with fixity 2, and this gives the possibilities
in Lines 5c and 5d.
Similarly, if $q \equiv -1$ modulo $4$, then we obtain the branching data in Lines 5e and 5f.

The next case is $G=\Sz(q)$, where $q$ is a power of $2$ with odd exponent. There are two possibilities for a fixity 4 action, with point stabilisers of different orders, which means that there can be up to two orbits on which $G$ acts with fixity 4.
This gives Lines 8a-8c of table \ref{DataList}.
It is also possible for $G$ to act with fixity 2, but only with point stabilisers of order $q-1$, and this can only happen on up to two orbits.
Combining up to two orbits with fixity 4 and up to two orbits with fixity 2 gives Lines 8d-8i.

There is only one group left that can act with mixed fixity, namely $\M_{22}$. A unique non-regular orbit leads to Line 12a of Table \ref{DataList}, and it is also possible to have one orbit with fixity 4 and another one with fixity 3, as in Line 12b.  

The specifications about $g_0$, if there are any, stem from the fact that $G$ acts as a group of automorphisms on $X$ and therefore the Hurwitz formula must be satisfied.

\medskip
For the second statement we suppose that $l$ is one of the lists in Table \ref{DataList} and that it is a Hurwitz datum.
Then we replace $l$ by a list $l_0$ that contains the smallest possible cogenus (as indicated before Table \ref{DataList} and in the table) and otherwise has the same entries as $l$. We will prove the statement for $l_0$, and then it follows with Lemma 3.6 in \cite{SW} that the statement also holds for $l$.

We recall that the genus $g$ is at least 2 in all lists of the table. 
Using \texttt{GAP} (\cite{GAP4}) we calculate, for some individual groups and the corresponding lists in the table, a generating tuple that satisfies Lemma \ref{crit}. 
This is mostly based on information about the conjugacy classes of the groups and the function \texttt{GeneratingMCOrbits} from the \texttt{MapClass}
package (\cite{{JMSV18}}).
We only write a few of the generating tuples down, but all the code that we use is explained in a supplementary file. In all generic cases, we refer to one of our previous lemmas.

We begin with the case where $G=\Alt_7$. Then \texttt{GAP} gives us the generating tuple
$( ((1,2)(3,7,5)(4,6)), ((1,3,7,4)(2,6)), (3,4,5,6,7) )$,
which satisfies Lemma \ref{crit}, and therefore the list $[ \Alt_7, g, 1 \mid [5,1] ]$ has a witness.

For the list $[ \Alt_7, g, 1 \mid [5,1], [7,1] ]$ we find the generating tuple\\
$( ((1,3,6,4,2)), ((2,4)(6,7)), ((2,5,4,6,3), (1,2,3,4,5,7,6)) )$, again using \texttt{GAP}.\\[2ex]
Next we let $G:=\PSL_2(7)$. Here we have six different lists to look at, and 
in the supplementary file we explain how we find generating tuples for each of them, with \texttt{GAP}. We sometimes deviate slightly from our standard method, depending on the exact entries of the lists, because our generic function based on \texttt{GeneratingMCOrbits} might be very slow.

The next nine lines of the table belong to the group $G=\PSL_2(8)$, with different lists. A subtle detail needs to be taken into account for List 3c, but otherwise we can calculate suitable generating tuples very easily with \texttt{GAP}. \\[1ex]

For $G=\PSU_4(3)$ we have two lists to check, and then we move to our first series of generic examples, namely $\PSL_2(q)$. We note that $q$ is odd and $q \ge 9$, with different variations of branching data. 
For the lists 5a and 5b we turn to Lemma \ref{PSL21} if $q \notin \{9,11,19\}$, and we use Lemmas \ref{PSL22}, \ref{PSL23} and \ref{PSL24} for the cases where $q \in \{9,11,19\}$, respectively. 
Then, for the mixed fixity cases in the lines 5c--5f, we consider Lemmas  
\ref{PSL22} and \ref{PSL26}.

The next generic examples come from the series $\PSp_4(q)$ and $\POm^-_8(q)$, in Lines 6 and 7. Here we refer to Lemmas \ref{Sp4} and \ref{POmega}.

The series of Suzuki groups gives many different examples, and the lists in Lines 8a -- 8i are all covered by Lemma \ref{suz}.

For the series $\,^3D_4(q)$ and $\,^2G_2(q)$ we have one type of list, each, and these have witnesses by Lemmas \ref{d4} and \ref{G2}.

For $\M_{11}$, we have one list, 
namely $[ \M_{11}, g, 1 \mid [5,1] ]$, 
and \texttt{GAP} give us a suitable generating tuple: \\[1ex] 
$( (( 1, 4,11,10, 7, 5, 8, 6, 9, 3, 2)), (( 3, 7, 4, 5)( 6,11, 9, 8)),$\\
$(( 2, 5, 4, 7, 6)( 3,11, 9, 8,10)) )$.

\vspace{1ex}
There are two lists for $\M_{22}$ and one list for $\J_1$, and again we do not display the generating tuples that we found with \texttt{GAP}.
The group $\J_1$ is the only one that we treat differently from the others, by using the \texttt{TomLib} package (\cite{MNP24}) and constructing the group from there. 
This concludes the proof of Theorem \ref{main4fpH}.



\section{The \texttt{GAP} code}

All the \texttt{GAP} code that we use is based on the packages \texttt{MapClass} and, for $\J_1$, \texttt{TomLib}. We need three basic functions for our calculations:

\begin{lstlisting}[language=Gap, caption=Function \texttt{TestTuple}]
TestTuple:=function(G,cogen,tuple)
	local len,prod,prodcomm,i;
	
	# Product of the non-commutators
	len := 2*cogen;;
	if Length(tuple)-len > 0 then
		prod := Product(tuple{[len+1..Length(tuple)]});;
	else
		prod := One(G);;
	fi;
	
	#Product of the commutators
	if cogen=0 then
		prodcomm := One(G);;
	else
		prodcomm := One(G);;
		for i in [1..cogen] do
			prodcomm := prodcomm * Comm( tuple[2*i-1], tuple[2*i] );;
		od;
	fi;
	
	# Final test
	if prodcomm*prod = One(G) and Group(tuple)=G then
		return true;;
	else
		return false;;
	fi;
end;
\end{lstlisting}

The first function \texttt{TestTuple} checks whether a calculated tuple \texttt{tuple} satisfies Lemma \ref{crit} (b) and (c) for the group \texttt{G} and cogenus \texttt{cogen}. It should be mentioned that in the \texttt{GAP} code we use a different notation for tuples for our calculations. 

If $l=[G, g, g_0 \mid [m_1,n_1],\ldots, [m_r,n_r]]$ is a Hurwitz datum and $((a_1,\ldots, a_{g_0}),$\linebreak
$(b_1,\ldots, b_{g_0}), (c_{1,1}, \ldots c_{1,n_1},\ldots,c_{r,1},\ldots, c_{r,n_r}))$ is a generating tuple, then the notation for the generating tuple in \texttt{GAP} is \texttt{[ $a_1, b_1,\ldots, a_{g_0}, b_{g_0}, c_{1,1}, \ldots, c_{1,n_1}, \ldots,$}
\linebreak\texttt{$c_{r,1}, \ldots, c_{r,n_r}$ ]}.

The next function, named \texttt{GenTuple}, computes a generating set for a given Hurwitz datum. It uses the \texttt{GeneratingMCOrbits} function of the \texttt{MapClass} package, which means that this package must be loaded before using the function \texttt{GenTuple}.

\begin{lstlisting}[language=Gap, caption=Function \texttt{GenTuple}]
GenTuple:=function(G,cogen,orders)
	local ct,cc,list,pos,ro,tuple,orb,orb1,h,i,len;

	ct := CharacterTable(G);;
	cc := ConjugacyClasses(ct);;
	list := OrdersClassRepresentatives(ct);;
	pos := List(orders, i-> Positions(list,i) );;
	ro := List(pos, i-> Random(i) );;
	tuple := List(ro,i->Representative(cc[i]) );;

	orb := GeneratingMCOrbits(G,cogen,tuple : OutputStyle:="single-line");;
	len := Length(orb);;
	if len=0 then return; fi;

	for i in [1..len] do
		Print(Length(orb[i].TupleTable));
		Print(" tuples in orbit ");
		Print(i);
		Print("\n");
	od;

    Print("\n"); 
    Print("Picking random tuple in random orbit ...\n");
    orb1 := orb[Random([1..len])];;
    len := Length(orb1.TupleTable);;
    h := SelectTuple(orb1, Random([1..len]) );;

	Print("Testing tuple ... \n");
	if TestTuple(G,cogen,h) then
		Print("A generating tuple for the given Hurwitz datum is \n");
		return h;
	else
		Print("Something went wrong. Please Try again!\n");
		return;;
	fi;
end;
\end{lstlisting}

The function \texttt{GenTuple} expects three parameters, a group \texttt{G}, a cogenus \texttt{cogen} and a list \texttt{orders} of element orders of \texttt{G}. If $l=[G, g, g_0 \mid [m_1,n_1],\ldots, [m_r,n_r]]$ is a Hurwitz datum, then we run the function as

\begin{center}
    \texttt{GenTuple( $G$, $g_0$, [ $\underbrace{m_1,\ldots,m_1}_{n_1}, \ldots, \underbrace{m_r,\ldots,m_r}_{n_r}$])}.
\end{center}

The third function \texttt{GenTupleSilent} computes a generating tuple for a given Hurwitz datum in the same way as \texttt{GenTuple}, but the output of this function is reduced -- it only shows the generating tuple or an error message. The function \texttt{GenTupleSilent} expects the same parameters as \texttt{GenTuple}.

\begin{lstlisting}[language=Gap, caption=Function \texttt{GenTupleSilent}]
GenTupleSilent:=function(G,cogen,orders)
	local ct,cc,list,pos,ro,tuple,orb,orb1,h,i,len;

	ct := CharacterTable(G);;
	cc := ConjugacyClasses(ct);;
	list := OrdersClassRepresentatives(ct);;
	pos := List(orders, i-> Positions(list,i) );;
	ro := List(pos, i-> Random(i) );;
	tuple := List(ro,i->Representative(cc[i]) );;

	orb:=GeneratingMCOrbits(G,cogen,tuple : OutputStyle:="silent");;
	len := Length(orb);;
	if len=0 then return; fi;

	orb1 := orb[Random([1..len])];;
	len := Length(orb1.TupleTable);;
	h := SelectTuple(orb1, Random([1..len]) );;

	if TestTuple(G,cogen,h) then
		return h;;
	else
		Print("Something went wrong. Please try again!");
		return;;
	fi;
end;
\end{lstlisting}

Now we show the code for the specific groups and the lists from Table \ref{DataList}. 

To compute a generating tuple for the group $\Alt_7$ and for the Hurwitz datum 
$[\Alt_7, g, g_0 \mid [5, 1]]$, 
we run 
\verb|GenTuple(AlternatingGroup(7),1,[5]);|

For the list
and $[\Alt_7, g, g_0 \mid [5, 1], [7, 1]]$
we use a different approach, because \verb|GenTuple| is very slow. This is how we obtained the generating tuple in the proof of Theorem \ref{main4fpH}:

\begin{lstlisting}[language=Gap]
G := AlternatingGroup(7);;
tuple := GenTupleSilent(G,0,[5,5,7]);;
x := tuple[1];;
y := tuple[2];;
z := tuple[3];;
powerofx := List([1..Order(x)], i->x^i);;
ccxwox:=Filtered(G,i->Order(i)=Order(x) and  not (i in powerofx));;
l := Filtered(ccxwox, i -> i*x in ccxwox);;
a := Random(l);;
b := RepresentativeAction(G,a,a*x);;
tuplenew := [a,b,y,z];;


# Testing tuple
if TestTuple(G,1,tuplenew) then
	Print("A generating tuple for the given Hurwitz datum is:\n");
	Print(tuplenew);
else
	Print("Something went wrong. Please try again!\n");
fi;
\end{lstlisting}

For $\PSL_2(7)$ we have six different lists to check.
The ones that can be treated with the function 
\texttt{GenTuple} are 
$[\PSL_2(7), g, g_0 \mid [2, 1], [3, 1]]$,
$[\PSL_2(7), g, g_0 \mid [2, 1], [7, 1]]$ and 
$[ \PSL_2(7), g, g_0 \mid [2,1], [3,1], [7,1] ]$.
Bearing in mind the smallest value for $g_0$, 
we use, after setting $G$ to
$\PSL_2(7)$, the commands 

\texttt{GenTuple(PSL(2,7),1,[2,3]);},  
\texttt{GenTuple(PSL(2,7),1,[2,7]);} and\\ 
\texttt{GenTuple(PSL(2,7),0,[2,3,7]);} respectively.
  
A generating tuple for the list $[ \PSL_2(7), g, g_0 \mid [2,1], [3,1], [7,1] ]$ can be used to construct a tuple for the list 
$[ \PSL_2(7), g, g_0 \mid [2,1], [3,2], [7,1] ]$ as follows:

\begin{lstlisting}[language=Gap]
tuple := GenTuple(PSL(2,7),0,[2,3,7]);;
Print([tuple[1],tuple[2]^(-1),tuple[2]^(-1),tuple[3]]); 
\end{lstlisting}

Now we look at two cases that we treat differently.
The first list in the table is $[\PSL_2(7), g, g_0 \mid [2, 1]]$,
and we check it for the smallest possible value of $g_0$ as follows:

\begin{lstlisting}[language=Gap]
G := PSL(2,7);;
tuple := GenTupleSilent(G,1,[2,2]);;

a := Random(Filtered(G , i -> IsConjugate(G, i, i*tuple[3] )));;
b := RepresentativeAction( G, a, a*tuple[3] );;
tuplenew := [tuple[1],tuple[2],a,b,tuple[4]];;


# Testing tuple
if TestTuple(G,2,tuplenew) then
	Print("A generating tuple for the given Hurwitz datum is:\n");
	Print(tuplenew);
else
	Print("Something went wrong. Please try again!\n");
fi;
\end{lstlisting}

For the list $[\PSL_2(7), g, g_0 \mid [2, 1], [3, 2]]$,
we also do something different because GenTuple is very slow here:

\begin{lstlisting}[language=Gap]
G := PSL(2,7);;
tuple := GenTupleSilent(G,0,[2,2,3,3]);;

a := Random(Filtered(G , i -> IsConjugate(G, i, i*tuple[1] )));;
b := RepresentativeAction( G, a, a*tuple[1] );;
tuplenew := [a,b,tuple[2],tuple[3],tuple[4]];;

# Testing tuple
if TestTuple(G,1,tuplenew) then
	Print("A generating tuple for the given Hurwitz datum is:\n");
	Print(tuplenew);
else
	Print("Something went wrong. Please try again!\n");
fi;  
\end{lstlisting}

For $\PSL_2(8)$ we have even more cases to consider.
First we look at the lists that can be treated with \texttt{GenTuple}, again with the smallest cogenus, namely 
$[ \PSL_2(8), g, 1 \mid [2,1], [7,1] ]$,
$[ \PSL_2(8), g, 1 \mid [2,1], [9,1] ]$,\\
$[ \PSL_2(8), g, 0 \mid [2,1], [9,2] ]$ and 
$[ \PSL_2(8), g, 0 \mid [2,1], [7,1],[9,1] ]$.

Here we run the commands\\
\texttt{GenTuple(PSL(2,8),1,[2,7]);},
\texttt{GenTuple(PSL(2,8),1,[2,9]);},\\
\texttt{GenTuple(PSL(2,8),0,[2,9,9]);} and 
\texttt{GenTuple(PSL(2,8),0,[2,7,9]);}.

We can use the last generated tuple for some more lists as follows:

\begin{lstlisting}[language=Gap]
G := PSL(2,8);;
tuple := GenTupleSilent(G,0,[2,7,9]);
x := tuple[1];;
y := tuple[2];;
z := tuple[3];;
Print("7th Hurwitz datum: [ G, g, 0 | [2,1], [7,2], [9,1] ]\n");
Print("A generating tuple for the given Hurwitz datum is:\n");
Print([x,y^2,y^6,z]);

Print("8th Hurwitz datum: [ G, g, 0 | [2,1], [7,1], [9,2] ]\n");
Print("A generating tuple for the given Hurwitz datum is:\n");
Print([x,y,z^2,z^8]);

Print("9th Hurwitz datum: [ G, g, 0 | [2,1], [7,2], [9,2] ]\n");
Print("A generating tuple for the given Hurwitz datum is:\n");
Print([x,y^2,y^6,z^2,z^8]);   
\end{lstlisting}

The list $[ \PSL_2(8), g, 2 \mid [2,1] ]$ is treated differently:

\begin{lstlisting}[language=Gap]
G := PSL(2,8);;
tuple := GenTupleSilent(G,1,[2,2]);;

a := Random(Filtered(G , i -> IsConjugate(G, i, i*tuple[3])));;
b := RepresentativeAction(G,a,a*tuple[3]);;
tuplenew := [tuple[1],tuple[2],a,b,tuple[4]];;

# Testing tuple
if TestTuple(G,2,tuplenew) then
	Print("A generating tuple for the given Hurwitz datum is:\n");
	Print(tuplenew);
else
	Print("Something went wrong. Please try again!\n");
fi;
\end{lstlisting}

Finally, we construct a generating tuple for the list $[ \PSL_2(8), g, 0 \mid [2,1], [7,2] ]$. We must treat this differently from just running the function \texttt{GenTuple}, because in the definition of the variable \texttt{ro} in \texttt{GenTuple} we pick a random number. But for the list $[ \PSL_2(8), g, 0 \mid [2,1], [7,2] ]$, things go wrong if \texttt{ro[2]} and \texttt{ro[3]} are accidentally picked to be the same number within \texttt{GenTuple}. Therefore, we add an additional step to make sure that \texttt{ro[2]} and \texttt{ro[3]} are different.

\begin{lstlisting}[language=Gap]
G := PSL(2,8);;
cc := ConjugacyClasses(G);;
cco := List(cc, i-> Order(Representative(i)));;
pos := List([2,7,7], i-> Positions(cco,i));;
ro := List(pos, i->Random(i));;
if ro[2] = ro[3] then
	ro[3] := Random(Filtered(Positions(cco,7), i-> i <> ro[2]));;
fi;

reptup := List(ro, i -> Representative(cc[i]) );;
orb := GeneratingMCOrbits(G,0,reptup:OutputStyle:="single-line");;

orb1 := orb[Random([1..Length(orb)])];;
tuple := SelectTuple(orb1, Random([1..Length(orb1.TupleTable)]) );;

# Testing Tuple
if TestTuple(G,0,tuple) then
	Print("A generating tuple for the given Hurwitz datum is:\n");
	Print(tuple);
else
	Print("Something went wrong. Please try again!\n");
fi;
\end{lstlisting}

For the group $\PSU_4(3)$ we have two lists in our table.
Since \texttt{GenTuple} takes a long time here, we 
construct our generating tuples differently. First we look at the list $[ \PSU_4(3), g, 1 \mid [5,1] ]$, we calculate a candidate and check whether it satisfies all conditions.

\begin{lstlisting}[language=Gap]
G := PSU(4,3);;
x := First(G, i -> Order(i) = 5);;
a := First(G, i -> Order(i) = 12);;
cca := ConjugacyClass(G,a);;
a := First(cca, i-> i*x in cca and Group(i,x)=G);;
b := RepresentativeAction(G, a, a*x);;
tuple := [a,b,x^(-1)];;

# Testing Tuple
if TestTuple(G,1,tuple) then
	Print("A generating tuple for the given Hurwitz datum is:\n");
	Print(tuple);
else
	Print("Something went wrong. Please try again!\n");
fi;
\end{lstlisting}

For the second list, 
$[ \PSU_4(3), g, 1 \mid [5,1], [7,1] ]$,
we start similarly and then change the generators slightly:

\begin{lstlisting}[language=Gap]
G := PSU(4,3);;
x := First(G, i -> Order(i) = 5);;
a := First(G, i -> Order(i) = 12);;
cca := ConjugacyClass(G,a);;
a := First(cca, i-> i*x in cca and Group(i,x)=G);;	
b := RepresentativeAction(G, a, a*x);;

ccx := ConjugacyClass(G,x);;
y := First(ccx, i -> Order(i^(-1)*x^(-1))=7 );;
z := y^(-1) * x^(-1);;
tuple := [a,b,y,z];;

# Testing Tuple
if TestTuple(G,1,tuple) then
	Print("A generating tuple for the given Hurwitz datum is:\n");
	Print(tuple);
else
	Print("Something went wrong. Please try again!\n");
fi;
\end{lstlisting}

Finally, we have three sporadic groups to consider, with four lists in total.
\texttt{GenTuple} is often too slow here. For $\M_{11}$ and the list 
$[ \M_{11}, g, 1 \mid [5,1] ]$, we use the following code:

\begin{lstlisting}[language=Gap]
G := MathieuGroup(11);;
x := First(G, i -> Order(i) = 5);;
a := First(G, i -> Order(i) = 11);;
cca := ConjugacyClass(G,a);;
a := First(cca,i-> i*x in cca and Group(i,x)=G);;
b := RepresentativeAction(G, a, a*x);;
tuple := [a,b,x^(-1)];;

# Testing Tuple
if TestTuple(G,1,tuple) then
	Print("A generating tuple for the given Hurwitz datum is:\n");
	Print(tuple);
else
	Print("Something went wrong. Please try again!\n");
fi;    
\end{lstlisting}

The next lists in the table are for the group $\M_{22}$.
First we look at $[ \M_{22}, g, 1 \mid [5,1] ]$:

\begin{lstlisting}[language=Gap]
G := MathieuGroup(22);;
x := First(G, i -> Order(i) = 5);;
a := First(G, i -> Order(i) = 11);;
cca := ConjugacyClass(G,a);;
a := First(cca,i-> i*x in cca and Group(i,x)=G);;
b := RepresentativeAction(G,a,a*x);;
tuple := [a,b,x^(-1)];;

# Testing Tuple
if TestTuple(G,1,tuple) then
	Print("A generating tuple for the given Hurwitz datum is:\n");
	Print(tuple);
else
	Print("Something went wrong. Please try again!\n");
fi;
\end{lstlisting}

Next we consider $[ \M_{22}, g, 1 \mid [5,1],[7,1] ]$, again using some of the information from before.

\begin{lstlisting}[language=Gap]
G := MathieuGroup(22);;
x := First(G, i -> Order(i) = 5);;
a := First(G, i -> Order(i) = 11);;
cca := ConjugacyClass(G,a);;
a := First(cca,i-> i*x in cca and Group(i,x)=G);;
b := RepresentativeAction(G, a, a*x);;

ccx := ConjugacyClass(G,x);;
y := First(ccx, i -> Order(i^(-1)*x^(-1)) = 7 );;
z := y^(-1)*x^(-1);;
tuple:=[a,b,y,z];;

# Testing tuple
if TestTuple(G,1,tuple) then
	Print("A generating tuple for the given Hurwitz datum is:\n");
	Print(tuple);
else
	Print("Something went wrong. Please try again!\n");
fi;   
\end{lstlisting}

Finally, we have the code for the list 
$[\J_1, g, g_0 \mid [15, 1]]$. Here we use an additional \texttt{GAP} package, called \texttt{TomLib}. The \texttt{TomLib} package is only used for the construction of the group $\J_1$ as there is no predefined function in \texttt{GAP} to construct this group. After we generate $\J_1$ we proceed in a similar way as before.

\begin{lstlisting}[language=Gap]
LoadPackage("TomLib");;

tom := TableOfMarks("J1");;
G := UnderlyingGroup(tom);;
x := First(G, i -> Order(i) = 15);;
a := First(G, i -> Order(i) = 19);;
cca := ConjugacyClass(G,a);;
a := First(cca,i-> i*x in cca and Group(i,x)=G);;
b := RepresentativeAction(G, a, a*x);;
tuple := [a,b,x^(-1)];;

# Testing tuple
if TestTuple(G,1,tuple) then
	Print("A generating tuple for the given Hurwitz datum is:\n");
	Print(tuple);
else
	Print("Something went wrong. Please try again!\n");
fi;
\end{lstlisting}


	
	
\end{document}